\documentclass[12pt]{amsart}
\usepackage{amsmath,amscd,amssymb}
\usepackage{amscd}
\usepackage{colortbl}
\theoremstyle{plain}
\newtheorem{theorem}{Theorem}[section]
\newtheorem{lemma}[theorem]{Lemma}
\newtheorem{prop}[theorem]{Proposition}
\newtheorem{cor}[theorem]{Corollary}
\theoremstyle{remark}
\newtheorem{remark}[theorem]{Remark}
\newtheorem{examp}{Example}[section]
\theoremstyle{definition}
\newtheorem{definition}[theorem]{Definition}
\newcommand\lra{\longrightarrow}
\newcommand\PP[1][3]{\mathbb P^{#1}}
\newcommand\ZZ{\mathbb Z}
\newcommand\QQ{\mathbb Q}
\newcommand\FF{\mathbb F}
\newcommand\mo{\mathcal O}
\newcommand\mz{\mathbb Z}
\DeclareMathOperator{\Pic}{Pic}
\DeclareMathOperator{\codim}{codim}
\DeclareMathOperator{\mult}{mult}
\DeclareMathOperator{\tr}{tr}
\begin{document}

\title[Higher dimensional Calabi-Yau varieties]
{Construction and examples of higher-dimensional modular Calabi-Yau manifolds}
\author{S.~Cynk and K.~Hulek}
\subjclass[2000]{14G10, 14J32, 11G40}
\begin{abstract}
We construct several examples of higher-dimensional Calabi-Yau manifolds and prove their
modularity.
\end{abstract}
\maketitle

\section{Introduction}
\label{sec:intro}

As a consequence of Wiles' proof of the Taniyama-Shimura-Weil conjecture \cite{Wiles} there
has been considerable interest in the modularity of Calabi-Yau manifolds in recent years.

The case of dimension $2$ was first considered by Shioda and
Inose \cite{SI} who studied $K3$ surfaces with maximal Picard number, so-called singular
$K3$ surfaces. They showed that these surfaces can be defined over number fields and computed their
Hasse-Weil zeta-function. In the case of a singular $K3$ surface the transcendental lattice is
$2$-dimensional. If the surface is defined over $\QQ$, then Livn\'e \cite{Livne} showed 
that the corresponding $2$-dimensional Galois representation is related to a weight $3$
modular form. 

In dimension $3$ rigid Calabi-Yau manifolds are simplest in the sense that they have $2$-dimensional
middle cohomology. By a variant of the Fontaine-Mazur conjecture \cite{FM}, also asked by Yui
(see e.~g. \cite{Yui} for a recent account), one expects that the middle cohomology 
of a rigid Calabi-Yau threefold defined over $\QQ$, gives rise to an $L$-series, which is that of
a weight $4$ modular form. After numerous examples by various authors were exhibited,
Dieulefait and Manoharmayum \cite{DM} proved the modularity conjecture for rigid 
Calabi-Yau threefolds under mild conditions on the primes of bad reduction. Examples and results
about non-rigid modular Calabi-Yau threefolds can be found e.g. in \cite{HV1}, \cite{HV2}.
For a very recent survey, including lists of practically all known examples we refer the reader to 
the book by Meyer \cite{bookofMeyer}.

However, practically no examples seem to be known in higher dimension, and it is the aim of
this paper to fill this gap. The first type of examples we give, arises inductively from the
Kummer construction described in Proposition \ref{prop:inv1} in section \ref{sec:kummer}. 
The manifolds obtained in this way are resolutions of quotients of products of Calabi-Yau manifolds
by a group of the form $\ZZ_2^n$.
With this method one can construct 
several examples of modular Calabi-Yau manifolds (in any
dimension). The middle cohomology 
(if the dimension is odd), resp. the transcendental lattice (if the dimension is even) is a tensor 
product of the middle cohomologies of modular Calabi-Yau manifolds of lower dimension. In some 
cases this tensor product (or, more precisely, its semi-simplification) splits into $2$-dimensional  
modular pieces. 

In order to obtain higher dimensional Calabi-Yau manifolds with small (e.~g. $2$-dimensional) middle
cohomology, one has to refine the Kummer construction by taking quotients with respect to 
bigger groups. We consider suitable actions of the groups  $G=\ZZ_3^n$ or $\ZZ_4^n$ and discuss 
this in particular in the case of quotients of the form $(E \times \ldots \times E)/G$, where 
$E$ is an elliptic curve with extra automorphisms (see Sections \ref{sec:triple} and \ref{sec:foufold}). 
We show that these quotients have a smooth Calabi-Yau model, whose
middle cohomology (if the dimension is odd), resp. the transcendental lattice (if the dimension is even)
is $2$-dimensional. Moreover we show modularity and determine the corresponding cusp forms. 

The final example, which we discuss, goes back to Ahlgren \cite{ahlgren}. He considers a $5$-dimensional
affine variety $X$ which is a double cover of $5$-space branched along $12$ hyperplanes and relates
the number of points of $X(\FF_p)$ to the cusp form $g_6(q)=\eta(q^{2})^{12}$ of weight $6$ and level $4$. 
We prove in Theorem
\ref{theo:Ahlgren} that $X$ has a smooth projective model which is a $5$-dimensional Calabi-Yau manifold
with $b_1=b_3=0$ and $b_5=2$, whose $L$-series of the middle cohomology is 
that of the weight form $g_6$.

\subsection*{Acknowledgements}

{\small
We are grateful to the DFG for support under
grant Hu 337/5-2 in the frame of the Schwerpunktprogramm SPP 1084 ``Globale Methoden in der
komplexen Geometrie''. This grant supported the stay of the first named author at the
University of Hannover, who would like to thank this University for kind hospitality and
excellent working conditions. 
We would also like to thank M.~Sch\"utt for numerous discussions.
}

\section{The Kummer construction}
\label{sec:kummer}

We start by generalizing the Kummer construction, which has been used to construct Calabi-Yau
threefolds as quotients of the product of a $K3$ surface with an involution and an elliptic curve
modulo the diagonal involution.
To begin with, let $Y$ be a projective manifold of dimension $n$ with $H^{q}(\mo_{Y})=0$ for $q > 0$ and let
$D \in |-2K_{Y}|$ be a smooth divisor. The line bundle $-K_Y$ defines a double covering
\[
\pi: X \longrightarrow Y
\] 
branched along the divisor $D$, and
\[
K_X=\pi^{\ast}(K_Y + (-K_Y)) =0.
\]
Moreover, since
\[
\pi_{\ast}(\mo_{X})=\mo_{Y} \oplus K_Y
\]
it follows that for $0<q<n$:
\[
H^{q}(\mo_{X})\cong H^{q}(\mo_{Y})\oplus
H^{q}(K_{Y})\cong H^{q}(\mo_{Y})\oplus
H^{n-q}(\mo_{Y}) = 0
\]
and, therefore, the variety $X$ is a Calabi-Yau manifold.

Now assume that we have a pair $Y_{i}, i=1,2$ of
algebraic manifolds of
dimension  $n_i$, together with smooth divisors $D_{i}\in |-2K_{Y_{i}}|$. 
Moreover, assume that $H^{q}(\mo_{Y_{i}})=0$ for $i=1,2$ and $q>0$ and let $X_i$ be the double covers
described above. By construction, the product $X_1 \times X_2$ admits an action 
of $\ZZ_2 \oplus \ZZ_2$.

\begin{prop}\label{prop:inv1}
Under the above assumptions the quotient of the product $X_{1}\times X_{2}$ by the diagonal
involution admits a crepant resolution $X$, 
which is a (smooth) Calabi-Yau manifold. Moreover, there is a double cover $X \to Y$, branched
along a smooth divisor $D$ with $H^{q}(\mo_{Y})=0$ for $q>0$.  
\end{prop}

\begin{proof}
 The resolution may be described as follows: we start with the blow-up 
\[\sigma : Y\lra Y_{1}\times Y_{2}\]
of
$Y_{1}\times Y_{2}$ along $D_{1}\times D_{2}$. Denote the exceptional divisor by $E$ and
let 
\[
D=\sigma^{*}(D_{1}\times Y_{2}+Y_{1}\times D_{2})-2E
\]
be the strict transform of $D_{1}\times Y_{2} \cup Y_{1}\times D_{2}$. 
Since $D_{1}\times Y_{2}$ and $Y_{1}\times D_{2}$ intersect
transversally along $D_{1}\times D_{2}$, the divisor $D$ is smooth and 
iso\-mor\-phic to the disjoint union of $D_{1}\times Y_{2}$ and $Y_{1}\times
D_{2}$.

Moreover $D=\sigma^{*}(D_{1}\times Y_{2}+Y_{1}\times D_{2})-2E \sim
\sigma^{*}(\pi_{1}^{*}(-2K_{Y_{1}})+\pi_{2}^{*}(-2K_{Y_{2}}))-2E=
\sigma^{*}(-2K_{Y_{1}\times Y_{2}})-2E \sim -2K_{Y}$. Since $Y$ and
$Y_{1}\times Y_{2}$ are smooth birational projective manifolds 
$H^{q}(\mo_{Y})\cong H^{q}(\mo_{Y_{1}\times
Y_{2}})=0$ for $q>0$ (by the K\"unneth decomposition), and hence the double cover $X$ of $Y$ branched
along $D$ is a Calabi-Yau manifold.

Clearly, $X$ is birational to the quotient of $X_{1}\times X_{2}$ by
the action of the diagonal in $\mz_{2}\oplus\mz_{2}$. More precisely, the 
fixed point set of the diagonal involution is the inverse image $B$ of $D_{1}\times D_{2}$. 
The quotient by the diagonal involution has transversal $A_{1}$-singularities along the image of $B$. 
Let $Z$ be the blow-up of $X_{1}\times X_{2}$ along $B$. The involution 
lifts to $Z$ with fixed point set equal to the exceptional divisor
$\tilde B$,  which 
is ruled over $B$. The quotient of $Z$ by this 
involution is isomorphic to $X$, i.e., we have a commutative diagram of the form
\[
\begin{CD}
  X_1 \times X_2/ \mz_{2}  @<<< Z/ \mz_{2} =X\\
  @VVV @VV\pi V\\
Y_1 \times Y_2 @<<< Y
\end{CD}
\]
where the horizontal lines are inverse maps to blow-ups and the vertical lines are branched double covers.
\end{proof}

The above proposition allows us to use the covering $X \to Y$ inductively, and thus to
construct higher-dimensional Calabi-Yau manifolds.

The Euler characteristic of $X$ depends not only on the Euler
characteristics of $X_{1}$ and $X_{2}$, but also on the involution. By standard topological 
arguments we obtain
\[e(X)=\frac12 e(X_{1})e(X_{2})+\frac32e(D_{1})e(D_{2})\]
and
\[e(D)=\frac12e(X_{1})e(D_{2})+\frac12e(D_{1})e(X_{2})+e(D_{1})e(D_{2}).\]

In the special case, when $Y_{2}$ is an elliptic curve, branched over $4$ points in $\PP[1]$, we have
\[e(X)=6e(D_{1}),\quad e(D)=2e(X_{1})+4e(D_{1}).\]

\begin{examp}
The case where $X_1$ is a $K3$-surfaces with an involution whose quotient is rational and
$X_2$ is an elliptic curve was studied independently by Borcea \cite{Bor2} and Voisin \cite{Voisin}
in the context of mirror symmetry.
Already in this case we have several 
possibilities leading to different Euler numbers. Namely,
$e(D_{1})$ is an even integer ranging from $-18$ (for a smooth plane
sextic) to  $20$ ($10$ lines coming from the resolution of six lines in
$\PP[2]$ with four triple points). 
If  $X_{2}$ is an elliptic curve we get 
Calabi-Yau $3$-folds with Euler numbers equal to $-108$, $-96$, $-84$,
$-72$, $-60$, $-48$, $-36$, $-24$, $-12$, $0$, $12$, $24$, $36$, $48$, $60$, $72$, $84$, $96$, $108$,
$120$. 
  
\end{examp}

The Hodge numbers of $X$ cannot be computed in a similarly straightforward
way. If we know the Hodge numbers of   $X_{1}$ and $X_{2}$, we can
compute the Hodge numbers of $X_{1}\times X_{2}$. The involution will kill the
skew-symmetric part of the Hodge groups and preserve the symmetric
part. But we also have to take into account the contribution to the cohomology coming from 
the blow-up of $B$ and describe the action of the involution on it.

\begin{prop}\label{prop:inv2}
Let $X_{1},\dots,X_{n}$ be Calabi-Yau manifolds with involutions as above.
The quotient of the product   $X_{1}\times\dots\times X_{n}$ by the
action of $\{(m_1,\dots,m_n)\in \mz_{2}^{n}| 
m_{1}+\dots+m_{n}=0\}\cong\mz_{2}^{n-1}$ has a crepant resolution of
singularities which is a Calabi-Yau manifold.
\end{prop}
\begin{proof}
We shall proceed by induction on $n$. The case $n=2$ follows from
Proposition~\ref{prop:inv1}. Since the resulting Calabi-Yau manifold
has again an involution, we can iterate the procedure. For a sequence of
Calabi-Yau manifolds $X_{i}$ with involution we have the following
factorization 
\begin{eqnarray*}
(X_{1}\times
\dots\times X_{n})/\mz_{2}^{n-1}\cong
\left((X_{1}\times
\dots\times X_{n})/\mz_{2}^{n-2}\right)/(\mz_{2}^{n-1}/\mz_{2}^{n-2}) 
\end{eqnarray*}
 where $\mz_{2}^{n-2}$ denotes the group $\{(m_1,\dots,m_n)\in \mz_{2}^{n}| 
m_{1}+\dots+m_{n-1}=m_{n}=0\}$. 
Consequently 
\[(X_{1}\times \dots\times X_{n})/\mz_{2}^{n-1}\cong
\left(\left((X_{1}\times
\dots\times X_{n-1})/\mz_{2}^{n-2}\right)\times X_{n}\right)/\mz_{2},\] 
which proves the proposition.  
\end{proof}

\begin{cor}
Let $E_{i}; i=1,\dots,n$ be elliptic curves. The quotient
$E_{1}\times\dots\times E_{n}$ by the action of $\mz_{2}^{n-1}$ has a
smooth model $X^n$ which is a Calabi-Yau manifold with Euler characteristic
$e(X^{n})=\frac12(6^{n}+3(-2)^{n})$. 
\end{cor}

We would like to remark that quotients of the form $(E_{1}\times E_2 \times E_3)/ \mz_{2}^{2}$ were
first considered by Borcea \cite{Bor}, who also proved that the resulting Calabi-Yau threefolds have
CM if and only if the factors $E_i$ have CM. 

\begin{lemma} \label{invariantpart}
If $n$ is odd, then 
\[
H^{n}(X^{n})\cong
H^{n}(E_{1}\times\dots\times
E_{n})^{\mz_{2}^{n-1}}\cong H^{1}(E_{1})\otimes\dots\otimes
H^{1}(E_{n}).
\]
For $n$ even the (invariant) submotive
$H^{n}(E_{1}\times\dots\times
E_{n})^{\mz_{2}^{n-1}}$ of $H^{n}(X^{n})$ is isomorphic to the direct 
sum of a submotive generated by cycles of products of $n/2$ fibres 
and a submotive $I(X^n)\cong H^{1}(E_{1})\otimes\dots\otimes 
H^{1}(E_{n})$. The motive $I(X^n)$ contains the transcendental 
submotive, i.e., the orthogonal complement to the algebraic cycles of 
$X^n$.
\end{lemma}
\begin{proof}
We first consider the invariant part of the middle cohomology 
of $E_{1}\times\dots\times E_{n}$. Any tensor product $\otimes_{j} H^{i_j}(E_j)$ which contributes 
to this must have $\sum_j{i_j}=n$. Now assume that at least one $i_j=1$. 
Then we must have that all $i_j=1$, since
otherwise one can find some $\varepsilon \in {\mz_{2}^{n-1}}$ which
acts by  $-1$ on $\otimes_{j} H^{i_j}(E_j)$. If $n$ is odd, 
then $\sum_j{i_j}=n$ can only occur if at least one, and hence, by 
the above argument, all $i_j=1$. We finally remark that $X^n$ is of the form $Z_n/ {\mz_{2}^{n-1}}$
where $Z_n$ arises from the product $E_{1}\times\dots\times E_{n}$ by blowing-up rational 
submanifolds. This only contributes to the even cohomology and this contribution 
is spanned by algebraic cycles.
\end{proof}

This discussion easily implies the

\begin{prop}\label{prop:invmod}
Assume that the $E_i$ are defined over $\QQ$ with the involution given
as $x\mapsto -x$ and let $L(X^n,s)$, resp. $L(I(X^n),s)$ be the $L$-series 
associated to the Galois action on
$H^{n}(X^{n})$ for $n$ odd and the submotive $I(X^{n})$ for $n$ even. Then 
\[L(X^n,s)=L(g_{E_{1}}\otimes\dots\otimes g_{E_{n}},s),\]
\[ \mbox{ resp. } 
L(I(X^n),s)=L(g_{E_{1}}\otimes\dots\otimes g_{E_{n}},s)\]
where the $g_{E_{i}}$ are the cusp forms associated to $E_{i}$. 
\end{prop}
\begin{proof}
The only statement, which requires a proof, is that $X_n$ is defined over $\QQ$. But this is clear, since
the factors $E_i$, the involutions and the locus which is blown up are all defined over $\QQ$. 
\end{proof}

Here we consider $g_{E_{1}}\otimes\dots\otimes g_{E_{n}}$ as the tensor product of  Galois-modules. For the analytic properties of (some) tensor products see~\cite{KS}.

\begin{remark}
For a generic choice of elliptic curves $I(X^n)$ equals the 
transcendental submotive of $X^n$, whereas in special cases it may be 
striclty bigger. For instance, if the factors $E_{2i-1}$ and $E_{2i}$ 
($i=1,\dots,n/2$) are isogeneous, then $I(X^n)$ contains the product of the
graphs of isogenies. If, moreover, the $E_i$'s have complex 
multiplication then $I(X^n)$ contains also the product of graphs of 
complex multiplications. Note that this is in agreement with
the appearance of the factors $L(s-\tfrac{n}{2})$ and 
$L(\chi_{-d},s-\tfrac{n}{2})$ in the $L$-series given below.
\end{remark}

We now specialize the situation even further and assume
that all $E_{i}$ are isomorphic to  an elliptic curve $E$ with 
complex multiplication in $\QQ(\sqrt{-d})$.
If $n$ is odd, then 
\[L(X,s)=L(g_{n+1},s)^{\binom n0}L(g_{n-1},s-1)^{\binom  n1}
\dots L\left(g_{2},s-\tfrac{n-1}2\right)^{\binom n{(n-1)/2}}
\]
and if $n$ is even, then
\[L(I(X),s)=L(g_{n+1},s)^{\binom n0}L(g_{n-1},s-1)^{\binom  n1}
\dots L\left(g_{3},s-\tfrac{n-2}2\right)^{\binom n{(n-2)/2}}%\times
\]
\[
\times L(\chi_{-d},s-\tfrac{n}{2})^{\frac12{\binom  n{n/2}}}L(s-\tfrac{n}{2})^{\frac12{\binom  n{n/2}}}.
\]

Here   $\zeta(\chi_{-d},s)$ is the Dirichlet $L$-function defined by the character
associated to the number field $K=\QQ(\sqrt{-d})$, i.e
$\chi_{-d}(p)=(\frac{-d}{p})$ and
$g_k$ is the cusp form corresponding to the $(k-1)$st power of the
Gr\"os\-sen\-cha\-rak\-ter $\psi$ of the elliptic curve $E$ \cite{Ribet}. The cusp form $g_{k}$ has
weight $k$ and  complex multiplication in the same field
as $E$. The Fourier coefficient $a_{n}(g_{k})$ is given by the sum
of the values of the Gr\"ossencharakter $\psi^{k-1}$ at the ideals in the ring
$\mo_{K}$ of integers in $K$ of norm $n$,
relatively prime to the conductor of $E$. For a prime $p$ which is
inert in $\mo_{K}$, we get $a_{p}=0$,  because there is no ideal in $\mo_{K}$ with norm $p$. For
a split prime $p$ we have $p=\alpha_p\bar\alpha_p$ for some $\alpha_p\in\mo_{K}$, which is determined by $E$.
Then
$a_{p}(g_{k})=\alpha_p^{k-1}+\bar\alpha_p^{k-1}$, 
more explicitly,
we have $a_{p}(g_{3})=a_{p}^{2}-2p$, $a_{p}(g_{4})=a_{p}^{3}-3pa_{p}$,
$a_{p}(g_{5})=a_{p}^{4}-4pa_{p}^{2}+2p^{2}$ and so on.

In terms of the associated Galois representations,
the connection between the forms $g_k$ and $g_2$ can be described as
follows. Consider the representation associated to $g_2$ and let
$(\alpha_p,\bar{\alpha}_p)$ be the eigenvalues of $\operatorname{Frob}_p$ for
primes $p$ with $\chi_{-d}(p)=1$. If $\chi_{-d}(p)=-1$, then the corresponding
eigenvalues are $(ip^{\frac{1}{2}},-ip^{\frac{1}{2}})$.
The eigenvalues of $g_k$ are then $({\alpha_p}^{k-1},{\bar{\alpha}_p}^{k-1})$ for $\chi_{-d}(p)=1$ and
$(p^{\frac{k-1}{2}},-p^{\frac{k-1}{2}})$ for $k$ odd and $\chi_{-d}(p)=-1$, resp.
$(ip^{\frac{k-1}{2}},-ip^{\frac{k-1}{2}})$ for $k$ even and $\chi_{-d}(p)=-1$.

We want to conclude this section by discussing one further example of our
Kummer construction. As the first factor we choose the rigid
Calabi-Yau 3-fold $X_{3}$, constructed as a resolution of
singularities of the double covering of $\PP$ branched along the following
arrangement of eight planes
\[xt(x-z-t)(x-z+t)y(y+z-t)(y+z+t)(y+2z)=0.\]
For a discussion of the properties of this (and other) double octics
see \cite{CM} and \cite[Octic Arr. No. $19$]{bookofMeyer}. 
As the second factor we take the $K3$ surface $S$ which is
obtained as a desingularization of the double sextic 
branched along the following arrangement of six
lines
\[xy(x+y+z)(x+y-z)(x-y+z)(x-y-z)=0.\]
Both $X_3$ and $S$ come with natural involutions which allow us to apply 
Proposition \ref{prop:inv1}. In this way we obtain
a smooth Calabi-Yau fivefold
$X_{5}$, which is the quotient of a blow-up $\widetilde{X_{3}\times S}$ of $X_{3}\times S$ by an
involution. So the Hodge groups of $X_{5}$ are the invariant part of
the Hodge groups of $\widetilde{X_{3}\times S}$ . Since we blow up products of lines and
blown-up planes, the odd-dimensional cohomology groups of  $\widetilde{X_{3}\times S}$ and
${X_{3}\times S}$ are the same.

Now, the odd-dimensional cohomology groups of $S$ vanish, whereas the only
odd-dimensional cohomology of $X_{3}$ is $H^{3}(X_{3})=H^{3,0}\oplus H^{0,3}$,
which is anti--invariant. The anti--invariant part of the co\-ho\-mo\-lo\-gy of $S$ is
$H^{2,0}\oplus H^{0,2}\cong T(S)\otimes_{\mathbb Z}\mathbb C$, where $T(S)$ is the transcendental
lattice.
Consequently, $b_{1}(X_{5})=b_{3}(X_{5})=0$ and $b_{5}(X_4)=4$ and moreover
\[H^{5}(X_{5})\cong H^{3}(X_{3})\otimes T(S).\]
Recall that (see  \cite{AOP} and \cite[p. 57]{bookofMeyer})
\[
L(T(S),s) \circeq L(g_3,s), \quad L(X_3,s) \circeq L(g_4,s)
\]
where $g_3$ and $g_4$ are the unique weight $3$, resp. weight $4$ Hecke eigenforms of level $16$ and $32$ 
with complex multiplication by $i$.
As usual $\circeq$ denotes  
equality upto a finite number of Euler factors. In concrete terms
\[
g_3(q)= \eta(q^{4})^6 = q-6q^5+9q^9+10q^{13}-30q^{17}+ \ldots
\]
and
\[
\quad g_4(q)=q + 22q^5 - 27q^9 - 18q^{13} - 94q^{17} + 359q^{25} +\ldots
\]
where $\eta(q)=q^{\frac{1}{24}}\prod_{n=1}^{\infty}(1-q^n)$ is the Dedekind $\eta$-function.
Both of these forms can be derived from the unique weight $2$ level $32$ newform
\[
g_2(q)=\eta(q^{8})^2\eta(q^{4})^2=q - 2q^5 - 3q^9 + 6q^{13} + 2q^{17} + \ldots 
\]
by taking the second, resp. third power of the  Gr\"os\-sen\-cha\-rak\-ter of $\mathbb Q[i]$ 
given as $\psi((\alpha))=\alpha$ for $\alpha\in\mathbb Z[i],\quad \alpha\equiv1\mod 2+2i$.

Hence we obtain that the L-series
of $X_{5}$ is the product of the L-series associated to $X_{3}$ and
$S$ and we also find that it factors  as 
\[L(X_{5},s) \circeq L(g_{4}\otimes g_{3},s) \circeq L(g_6,s)L(g_2,s-2)\]
where $g_2$ is as above and $g_{6}$ is a level $32$ cusp form of weight $6$, namely
\[
g_6(q)=q - 82q^5 - 243q^9 - 1194q^{13} + 2242q^{17} + 3599q^{25} + \ldots
\]
which can be derived from $g_2$ by taking the fifth power of the Gr\"os\-sen\-cha\-rak\-ter.
Obviously, we can iterate this procedure to obtain modular Calabi-Yau manifolds of higher dimension
(with increasingly complex middle cohomology).

\section{Calabi-Yau manifolds with an endomorphism of order $3$}
\label{sec:triple}

We shall construct for any positive integer $n$ a Calabi-Yau
$n$-fold $X_{n}$ with an endomorphism of order 3 such
that $\dim H^{n}(X_{n})=2$ for $n$ odd and $\dim T(X_{n})=2$ for
$n$ even, where $T(X_{n})\subset H^{n}(X_{n})$ is the transcendental part. 
Moreover, we
shall show that the (semi-simplifications of) the Galois representation
on $H^{n}(X_{n})$ (resp. $T(X_{n})$) and the Galois
representation associated to a suitable cusp form with CM by $\sqrt{-3}$ are isomorphic.

Fix the primitive third root of unity $\zeta=e^{2\pi i/3}$.
Let $X_{1}$ and $X_{2}$ be two Calabi-Yau manifolds admitting
$\mz_3$-actions, which do not preserve the canonical form. Moreover, assume that
the fixed point set of the action on $X_{1}$ is a smooth divisor,
whereas on $X_{2}$ it is a disjoint union of a smooth divisor and a
smooth codimension  two submanifold. Fix an automorphism $\eta_{1}$ of
$X_{1}$ such that $\eta_{1}^{*}\omega _{X_{1}}=\zeta\omega_{X_{1}}$
and an automorphism $\eta_{2}$ of
$X_{2}$ such that $\eta_{2}^{*}\omega _{X_{2}}=\zeta^{2}
\omega_{X_{2}}$ such that they act on $X_1$ and $X_2$ as described above. 
Then $\eta_{1}$ is given locally near the branch-divisor on $X_{1}$
 as $(\zeta,1,1\dots)$, whereas 
$\eta_{2}$ is given locally  either as $(\zeta^{2},1,1\dots)$ near the branch divisor on
$X_{2}$ or as $(\zeta,\zeta,1\dots)$ near the codimension $2$ fixed locus.

On $X_{1}\times X_{2}$ we have an action of $\mz_3\oplus\mz_3$, and we
consider the action of $\mz_3$ on $X_{1}\times X_{2}$,
given by the automorphism $\eta=\eta_{1}\times \eta_{2}$. 

\begin{prop}\label{Z3action}
Under the above assumptions the quotient variety $X_{1}\times
X_{2}/\mz_3$ has a resolution of singularities $X$, which is a
Calabi-Yau manifold. The manifold $X$ admits a $\mz_3$-action which 
satisfies the same assumptions as for $X_2$.
\end{prop}
\begin{proof}
The singularities of  $X_{1}\times X_{2}/\mz_3$ correspond to the fixed locus of $\eta$,
which is the cartesian product of the fixed point sets of $\eta_{1}$ and
$\eta_{2}$. Consequently, we get two kind of singularities: a singular
codimension two stratum $W_1$, which is a transversal $A_{2}$-singularity, and a codimension
three stratum $W_2$, which is a transversal cone over a triple Veronese surface. Both types of
singularities admit a crepant resolution (described explicitly below), and we denote the resulting manifold
by $X$. Since the canonical form
on $X_{1}\times X_{2}$ is $\eta$-invariant it descends to the
quotient, and thus to the crepant resolution. Consequently, we get
$\omega_{X}\cong \mo_{X}$. 

Denote by $W_{1}$ (resp. $W_{2}$) the union of the codimension two (resp. three) strata
of the fixed point set of $\eta$ and
consider the blow-up $Z_{1}$ of $X_{1}\times X_{2}$ along $W_{1}\cup
W_{2}$. Then  $\eta$ lifts to $Z_{1}$ and the fixed
point set is a codimension two subvariety lying over $W_{1}$ and
a divisor over $W_{2}$. Let $Z_{2}$ be the blow-up of $Z_{1}$ along
the codimension two fixed submanifold. Again, the action of $\mz_3$
lifts to $Z_{2}$ and the fixed point set is a divisor. So the quotient
$Z$ of $Z_{2}$ by the action of $\mz _{3}$ is a smooth manifold and it
is a blow-up of $X$. (In terms of the $A_2$-singularity, the difference between
$Z$ and $X$ is, that we blow up the point of intersection of the two $(-2)$-curves 
which come from the resolution of the $A_2$-singularity.)
Now observe that 
\[
H^{0}(Z_2,\Omega^{q}_{Z_{2}})=H^{0}(X_1 \times X_2,\Omega^{q}_{X_1 \times X_2})=0
\] 
for $q\not=0,n_{1},n_{2},n_{1}+n_{2}$ and hence, by taking the invariant part
with respect to the action of $\eta$, we obtain that 
$H^{0}(Z,\Omega^{q}_{Z})=H^{0}(X,\Omega^{q}_{X})=0$
for $p\not=0,n_{1}+n_{2}$. This proves that $X$ is a (smooth) Calabi-Yau
manifold.

The action of $\mz_3\oplus\mz_3$ on $X_{1}\times X_{2}$ induces an
action of $\mz_3$ on  $X$ generated by the induced action of
$\operatorname{id}\times\eta_{2}$ to $X$. We shall study this action
in local coordinates. For the transversal $A_{2}$-singularity we can
find local coordinates on $X_{1}\times X_{2}$ in which the action is given as $(\zeta,
\zeta^{2})$. Note that for simplicity we shall omit the coordinates on
which $\eta$ acts trivially. The quotient map is given by
$(x_{1},x_{2})\mapsto(u_{1},u_{2},u_{3})=
(x_{1}^{3},x_{2}^{3},x_{1}x_{2})$, and the image has
equation $u_{3}^{3}=u_{1}u_{2}$. The resolution of singularities is given
by blowing up the submanifold $u_{1}=u_{2}=u_{3}=0$. In suitable charts on the blown up
surface, the quotient map is then given by
\begin{eqnarray*}
  (y_{1},y_{2})&\mapsto&(y_{1},y_{1}^{2}y_{2}^{3},y_{1}y_{2})\\
  (y_{1},y_{2})&\mapsto&(y_{1}^{3}y_{2}^{2},y_{2},y_{1}y_{2})\\
  (y_{1},y_{2})&\mapsto&(y_{1}^{2}y_{2},y_{1}y_{2}^{2},y_{1}y_{2}).
\end{eqnarray*}
The map from  $X_{1}\times X_{2}$ to the resolution of the quotient
in local analytic terms is given by
$
(x_{1}^{3},\frac{x_{2}}{x_{1}^{2}}), 
(\frac{x_{1}}{x_{2}^{2}},x_{2}^{3})$, 
$(\frac{x_{1}^{2}}{x_{2}},\frac{x_{2}^{2}}{x_{1}})$, depending on the charts we work in. 
The action of $ \operatorname{id}\times\eta_{2}$ is given on
$X_{1}\times X_{2}$ as $(1,\zeta^{2})$, so it lifts to $X$ as
$(1,\zeta^{2}), (\zeta^{2},1)$ or $(\zeta,\zeta)$ respectively, depending
on the affine chart we consider. 

For the cone over the Veronese triple embedding the resolution is given by
the so-called canonical resolution \cite[(16.10, p. 199)]{Ueno}.
If the action on $X_{1}\times X_{2}$ is given in local coordinates as
$(\zeta,\zeta,\zeta)$, then the map inverse to the resolution of the
quotient is given as
$\left(x_{1}^{3},\frac{x_{2}}{x_{1}},\frac{x_{3}}{x_{1}}\right)$, 
$\left(\frac{x_{1}}{x_{2}},x_{2}^{3},\frac{x_{3}}{x_{2}}\right)$ or 
$\left(\frac{x_{1}}{x_{3}},\frac{x_{2}}{x_{3}},x_{3}^{3}\right)
$.
The  action of $ \operatorname{id}\times\eta_{2}$ is given on
$X_{1}\times X_{2}$ as $(1,\zeta,\zeta)$, so it  lifts to $X$ as
$(1,\zeta,\zeta), (\zeta^{2},1,1)$ and
$(\zeta^{2},1,1)$ respectively. 

In all cases  $X$ satisfies the assumptions made for $X_{2}$.
\end{proof}

\begin{remark} We can now use the Calabi-Yau manifold $X$ together with the 
$\mz_3$-action on it to repeat this constructions inductively.
\end{remark}
 
We consider an elliptic curve defined over $\QQ$ with an automorphism
of order $3$, which we can, without loss of generality, assume to be in
Weierstrass form $y^{2}=x^{3}-D$. The automorphism $\eta$ is given by
$x\mapsto\zeta x$. 

\begin{theorem}\label{thm:level27}
Let $E$ be the elliptic curve with an automorphism $\eta$ of order $3$ and 
let $\bar X_{n}$ be the quotient of $E^{n}$ by the action of
the group 
\[\{(\eta^{a_{1}}\times\dots\times \eta^{a_{n}})\in
\operatorname{End}(E^{n}): a_{1}+\dots+a_{n}\equiv0\mod3\}.
\]
Then $\bar{X_{n}}$ has a smooth model $X_{n}$, which is a Calabi-Yau
manifold and $\dim(H^n(X_{n}))=2$, if $n$ is odd, resp. $\dim(T(X_{n}))=2$, if $n$ is even,
where $T(X_n)$ is the transcendental part of the cohomology.

Moreover, $X_n$ is defined over $\QQ$ and $L(H^n(X_{n}),s)\circeq L(g_{n+1},s)$, resp.
$L(T(X_n))\circeq L(g_{n+1},s)$, where
$g_{n+1}$ is the weight $n+1$ cusp form with complex multiplication in
$\QQ(\sqrt{-3})$, associated to the $n$-th power of the Gr\"ossencharakter of $E$. 
\end{theorem}
\begin{proof}
The claim about $X_n$ being a Calabi-Yau manifold follows by repeated application of 
Proposition \ref{Z3action}. To compute the middle cohomology, resp. its transcendental part,
we first notice that it is enough to compute the invariant part of the 
cohomology of $E^n$. This follows, since the divisors which we
introduce by blowing up, are linear spaces blown-up in some subspaces,
so their cohomology is generated by algebraic
cycles. 
The subspace $\otimes H^{10}(E) \oplus \otimes H^{01}(E)$ is always invariant.
If $n$ is odd, then, by an argument 
similar to the one we used in the proof of Proposition \ref{invariantpart}  this is the 
only contribution to the invariant part of $H^n(E^n)$. If $n$ is even we have, in addition,
summands of the form $H^{i_1}(E) \otimes \ldots \otimes H^{i_n}(E)$, where $i_k= 0$ or $2$ and
$\sum i_k=n$, which are also generated by  algebraic cycles.
 
Now we turn to the arithmetic statements. We first note that $\bar{X_{n}}$ is defined over $\QQ$,
since it is defined over $\QQ(\sqrt{-3})$ and invariant under the
Galois group. Since we blow up in submanifolds defined over $\QQ$, the resolution $X_n$ is also
defined over $\QQ$.

The endomorphism $\eta$ induces endomorphisms
$\eta_{p}:E(\bar\FF_{p})\lra E(\bar\FF_{p})$ (for $p\not=3,
p \nmid D $),
also of order 3. The induced endomorphisms  $\eta_{p}$ 
have three fixed points, and hence the Lefschetz fixed point formula implies
$\tr\eta^{*}_{p}=-1$. Now, if $l\equiv 1\mod 6$, then $\QQ_{l}$
contains a primitive root of unity $\rho_{l}$ and the eigenvalues of
$\eta_{p}^{*}$ are powers of $\rho_{l}$ which sum up to $-1$, and are, therefore,
equal to $\rho_{l}$ and $\rho_{l}^{2}$. Denote by
$v_{1},v_{2}\in H^{1}_{\text{\'et}}(E_p)$ the corresponding
eigenvectors. 
It is easy to see that
the subspace of $H^{1}_{\text{\'et}}(E_{p})^{\otimes n}$ invariant
under the action of $\ZZ_{3}^{n}$ is generated by
$v_{1}^{\otimes n}=v_{1}\otimes\dots\otimes v_{1}$ and $v_{2}^{\otimes
n}=v_{2}\otimes\dots\otimes
v_{2}$, so we need to compute the images of the tensor power of Frobenius
on $v_{1}^{\otimes n}$ and $v_{2}^{\otimes n}$. 
To this end, we shall need to compute the action of Frobenius
$\operatorname{Frob}_{p}^{\ast}$ in the base $v_{1}, v_{2}$. 

We shall consider the cases $p\equiv 1,5 \mod 6$ separately. 
For $p\equiv 1 \mod 6$ the Frobenius map $\operatorname{Frob}_{p}^{\ast}$
commutes with $\eta_{p}^{\ast}$, so it acts as $v_{1}\mapsto \alpha_{p}v_{1}$
and $v_{2}\mapsto \bar\alpha_{p}v_{2}$, where $\alpha_{p}$ and
$\bar\alpha_{p}$ are the eigenvalues of
$\operatorname{Frob}_{p}^{*}$. Consequently, the eigenvalues of
Frobenius on the invariant part of
$H^{1}_{\text{\'et}}(E_p))^{\otimes n}$ equal $\alpha_{p}^{n}$ and
$\bar\alpha_{p}^{n}$. 

For $p\equiv 5 \mod 6$ we have
$\operatorname{Frob}_{p}^{*}\circ\eta_{p}^{*}=\left(\eta_{p}^{*}\right)^{-1}
\circ\operatorname{Frob}_{p}^{*}$, which easily implies that in the
base $v_{1},v_{2}$ Frobenius is given by the matrix
$\left(\begin{array}{cc}
0&\lambda\\-\frac p\lambda&0
\end{array}\right)$. Consequently, the action of Frobenius on the
invariant subspace of $H^{1}_{\text{\'et}}(E_p)^{\otimes n}$ equals 
$\left(\begin{array}{cc}
0&\lambda^{n}\\(-\frac p\lambda)^{n}&0
\end{array}\right)$, with eigenvalues equal to $\pm p^{n/2}$ for $n$ even,
and $\pm ip^{n/2}$ for $n$ odd.

Taking all the cases together, we see that the Galois representation on
$H^{n}(X^{n})$ for $n$ odd, resp. $T(X^{n})$ for $n$ even, has the same
eigenvalues as the representation associated to the cusp form
$g_{n+1}$ associated to the $n$-th power of the Gr\"ossencharakter of
the elliptic curve $E$.
\end{proof}

\begin{remark}
In the case where $E$ is given by the equation $y^2=x^3 - 1/4$ the form $g_2$ is the unique
weight $2$ newform of level $27$, namely
\[
g_2(q)=\eta(q^{9})^{2}\eta(q^{3})^{2}=q - 2q^4 - q^7 + 5q^{13} +4q^{16} - 7q^{19}+\dots
\]
In this case
\[
g_3(q)=q + 4q^4 - 13q^7 - q^{13} + 16q^{16} + 11q^{19} + 25q^{25} +\dots
\]
and
\[
g_4(q)=\eta(q^{3})^{8}=q-8q^{4}+20q^{7}-70q^{13}+64q^{16}+56q^{19}+\dots
\]
which are the unique level 27 and 9 forms of weight 3 and 4. 
Cusp forms $g_k$ correspond to powers of the 
Gr\"ossencharakter of the field $\mathbb Q(\sqrt{-3})$ given by
$\psi((\alpha))=\alpha$ for $\alpha\in\mathbb Z[\frac{1+\sqrt{-3}}2]$, $\alpha\equiv 1\mod 3$.
For other models of $E$ one obtains appropriate twists of 
these forms.
\end{remark}

\section{Calabi-Yau manifolds with an endomorphism of order $4$}
\label{sec:foufold}
In this section we shall construct a similar example as in the previous
section, but with an endomorphism of order $4$. 
Let $X_{1}$ and $X_{2}$ be two Calabi-Yau manifolds admitting
$\mz_{4}$-actions $\eta_{1}$ and $\eta_{2}$. Assume that the
fixed point set of $\eta_{1}$ is a divisor, near which the action has a
linearization of the form $(1,i)$, whereas the fixed point set of
$\eta_{2}$ is a disjoint union of submanifolds of codimension one, two
or three, near which the action has a linearization as $(-i,1,\dots)$,
$(-1,i,1,\dots)$ and $(i,i,i,1\dots)$ respectively. 

On $X_{1}\times X_{2}$ we have an action of $\mz_{4}\oplus\mz_{4}$, and we
consider the action of $\mz_{4}$ on $X_{1}\times X_{2}$
given by the automorphism $\eta=\eta_{1}\times \eta_{2}$. 

\begin{prop}\label{Z4action}
Under the above assumptions the quotient $X_{1}\times
X_{2}/\mz_{4}$ has a resolution of singularities $X$, which is a
Calabi-Yau manifold. The manifold $X$ admits a $\mz_{4}$-action which 
satisfies the same assumptions as for $X_2$.
\end{prop}

\begin{proof}
We shall show that the quotient admits a crepant resolution of
singularities. We shall consider
the three cases, depending on the codimension of the
component of the fix-point set of $\eta_{2}$, separately. 

Near the fixed divisor of $\eta_{2}$ the action of $\eta$ on $X_{1}\times
X_{2}$ is locally given by $(i,-i)$ (as in similar proofs before, we omit the
variables on which $\eta$ acts trivially). Consequently the quotient
is a transversal $A_{3}$ singularity along the singular subvariety,
which can be resolved by blowing-up twice.  In  local
coordinates,  the map from  $X_{1}\times X_{2}$ inverse to the  resolution
is given in affine charts as
$
\left(x^{4},\frac y{x^{3}}\right),
\left(y^{4},\frac x{y^{3}}\right),
\left(\frac{x^{3}}y,\frac {y^{2}}{x^{2}}\right),
\left(\frac{y^{3}}x,\frac {x^{2}}{y^{2}}\right)$
or
$
\left(\frac{x^{2}}{y^{2}},xy\right).
$
The action of $\operatorname{id}\times\eta_{2}$ on $X_{1}\times X_{2}$
has a linearization $(1,-i)$, so it lifts to the resolution  as  
$(1,-i)$, $(1,-i)$, $(i,-1)$, $(i,-1)$ or
$(-1,-i)$. In all the cases, except the last one, the action is
exactly as we assume for $X_{2}$, but in the last case the fixed point of
the action is $(0,0)$, which does not belong to the domain of the map.

Now consider the singularity corresponding to a codimension two fixed
stratum of $\eta_{2}$. Then the action on $X_{1}\times X_{2}$ has a local
linearization of the form $(i,i,-1)$. We first divide by
the square of $\eta$, which is an involution with fixed point set
of codimension two resulting in transversal $A_{1}$-singularities. These we resolve by
blowing-up the singular locus. The action of $\mz_{4}$ lifts to this
resolution again as an involution with a codimension two fixed
point set, leading once more to transversal $A_{1}$-singularities, which we resolve with a single
blow-up. Simple computations show, that in terms of local
coordinates, the map from  $X_{1}\times X_{2}$  inverse to the
resolution looks in local coordinates like
$
\left(x^{4},\frac z{x^{2}},\frac yx\right)$,
$
\left(z^{2},\frac {x^{2}}z,\frac yx\right)$,
$
\left(\frac{x^{2}}z,x^{2}z,\frac yx\right)$,
$
\left(\frac xy,x^{2}y^{2},\frac z{xy}\right)$,
$
\left(\frac xy,z^{2},\frac{xy}z\right)$,
$
\left(\frac xy,\frac{xy}z,xyz\right)$,
$
\left(y^{4},\frac z{y^{2}},\frac xy\right)$,
$
\left(z^{2},\frac{y^{2}}z,\frac xy\right)$ or
$
\left(\frac{y^{2}}z,y^{2}z,\frac xy\right)
$.
The action of $\operatorname{id}\times\eta_{2}$ on $X_{1}\times X_{2}$
is linearized by $(1,i,-1)$, so it lifts to the resolution  as  
$(1,-1,i)$, $(1,-1,i)$, $(-1,-1,i)$, $(-i,-1,i)$,
$(-i,1,-i)$, $(-i,-i,-i)$, $(1,1,-i)$, $(1,1,-i)$ or $(1,1,-i)$.
The lifting satisfies the assumption made for $X_{2}$ in all
except the cases $3, 4, 5$ and $6$, when the fixed points do not lie in the
domain of the map. 

The last case is the fixed point stratum of $\eta_{2}$ of codimension $3$,
so the action on $X_{1}\times X_{2}$ has a local linearization of the
form $(i,i,i,i)$. Here again, it is easier to resolve in one step. 
On the quotient we get a transversal cone over the
Veronese fourfold embeding of $\PP[3]$. The crepant resolution is
given by the so-called canonical resolution \cite[(16.16,
p. 199)]{Ueno} for which the inverse map is  given as 
$\left(x^{4},\frac yx,\frac zx, \frac tx\right)$,
$\left(\frac xy,y^{4},\frac zy, \frac ty\right)$,
$\left(\frac xz,\frac yz,z^{4}, \frac tz\right)$ and
$\left(\frac xt,\frac yt, \frac zt,t^{4}\right)$, 
and so the action of $\operatorname{id}\times\eta_{2}$ lifts as
$(1,i,i,i)$, $(-i,1,1,1)$, $(-i,1,1,1)$ and $(-i,1,1,1)$ respectively, which
completes the proof.
\end{proof}

To produce an explicit example, consider the elliptic curve
$E$ given by the Weierstrass equation $y^{2}=x^{3}-Dx$,
where $D$ is a square-free integer. This curve has  complex
multiplication in the field $\QQ[i]$ and the map $\rho:(x,y)\mapsto (-x,iy)$
is an endomorphism of $E$ of order $4$. 

\begin{theorem}\label{thm:level32}
Let  ${\bar X}_{n}$ be the quotient of $E^{n}$ by the action of
the group 
\[\{(\eta^{a_{1}}\times\dots\times \eta^{a_{n}})\in
\operatorname{End}(E^{n}): a_{1}+\dots+a_{n}\equiv0\mod4\}.
\]
Then ${\bar X}_{n}$ has a smooth model $X_{n}$, which is a Calabi-Yau
manifold and $\dim(H^n(X_{n}))=2$ if $n$ is odd, resp. $\dim(T(X_{n}))=2$ if $n$ is even,
where $T(X_n)$ is the transcendental part of the cohomology.

Moreover,  $X_n$ is defined over $\QQ$ and $L(H^n(X_{n}),s)\circeq L(g_{n+1},s)$, resp.
$L(T(X_n))\circeq L(g_{n+1},s)$, where
$g_{n+1}$ is a weight $n+1$ cusp form with complex multiplication in $\QQ(i)$. 
\end{theorem}
\begin{proof}
  The existence of a crepant resolution follows from repeated application of Proposition~\ref{Z4action}, and the
  remaining statements can be proved exactly in the same way as in the proof of
  Theorem~\ref{thm:level27}. 
\end{proof}

\section{The example of Ahlgren}
\label{sec:example1}

Let $\bar X$ be the double cover of $\mathbb P^{5}$
branched along the union of the twelve hyperplanes
\[x(x-u)(x-v)y(y-u)(y-v)z(z-u)(z-v)t(t-u)(t-v)=0.\]
This is a projective closure of the fivefold studied by Ahlgren
\cite{ahlgren}. He proved that the number of
points defined over $\mathbb F_{p}$ on the affine part ($u=1$) of
this variety equals 
\[N(p)=p^{5}+2p^{3}-4p^{2}-9p-1-a_{p},\] 
where $a_{p}$ is the $p$-th Fourier coefficient of the unique normalized weight $6$ and
level $4$ cusp form (which is equal to $\eta^{12}(q^2)$).

Our goal here is to prove the following 

\begin{theorem} \label{theo:Ahlgren}
The variety $\bar X$ has a smooth model $X$ (defined over $\QQ$), 
which is a Calabi-Yau fivefold with Betti numbers 
  $b_{1}(X)=b_{3}(X)=0$, $b_{5}(X)=2$. More precisely, 
  $h^{50}=h^{05}=1$, $h^{14}=h^{23}=h^{32}=h^{41}=0$. The
  (semi-simplifications of the) Galois representation of the action of
  Frobenius on $H^{5}(X)$ and the Galois representation
  corresponding to the unique normalized cusp form of level $4$ and weight $6$ (which is $\eta^{12}(q^2)$),
  are isomorphic.
\end{theorem}

Before we can give the proof we need some preparations.
  
The variety $\bar X$ is a double cover of a degree twelve arrangement,
in the sense of Definition~\ref{def:arr} (see 
subsection~\ref{sec:arr} at the end of this section, where we collect the necessary
statements). In
Proposition~\ref{prop:arr} we describe a procedure, how to resolve
singularities of such a double cover, and our goal here is to check that
the arrangement satisfies the assumptions of that proposition. 

We shall distinguish the singularities by their multiplicity and dimension and denote the
resulting classes by $T_k$. Let $N_{k}$ be the number of
singularities of type $T_{k}$ that contain a given singularity. 
Then the situation can be summed up by table 1.
\begin{table}
\[
\begin{array}{|c|c|c|c|c|c|c|c|c|c|}\hline
\text{Type}&\text{dim}&\text{mult}&\text{\#}&N_{1}&N_{2}&N_{3}&N_{4}&N_{5}&N_{6}\\
\hline T_{1}&3&2&66&0&0&0&0&0&0\\
\hline T_{2}&2&3&148&3&0&0&0&0&0\\
\hline T_{3}&2&4&18&6&0&0&0&0&0\\
\hline T_{4}&1&4&117& 6&4&0&0&0&0\\
\hline T_{5}&1&5&36&10& 6& 1&0&0&0\\
\hline T_{6}&1&6&18&15&8&3&0&0&0\\
\hline T_{7}&0&5&12&10&10&0&5&0&0\\
\hline T_{8}&0&6&18&15&16&1&6&2&0\\
\hline T_{9}&0&7&12&21&23&3&8&3&1\\
\hline T_{10}&0&8&3&28&32&6&16&0&4\\
\hline T_{11}&0&9&4&36&21&9&9&9&6\\
\hline
\end{array}
\]
\caption{}
\end{table}
We see that $T_{2}, T_{4}, T_{5}, T_{7}, T_{8}, T_{9}$ are near pencil,
whereas $T_{1}, T_{3}, T_{6}, T_{10}$ and $T_{11}$ satisfy 
$\left\lfloor\frac {m(C)}2\right\rfloor=n-d(C)-1$. Hence $\bar X$ has a crepant resolution
of singularities $X$, which is 
a smooth Calabi-Yau variety. 

Studying the
singularities in the above table, we see that the only
prime of bad reduction is $2$ (due to taking the double cover). 
The exterior powers of the matrix of coefficients of the arrangement of
hyperplanes have coefficients equal to $0,\pm1$, so the reduction modulo an 
odd prime has the same number and type of singularities as in
characteristic $0$. Consequently, the same blow-ups as in
characteristic $0$ give a resolution of singularities. 

To prove modularity of $X$, we have to study the number of points of
$X_{p}$ in $\mathbb F_{p}$. In principle, it should be possible to 
give an explicit formula, as was done in the analogous situation in dimension $3$. However, in this case
there are many more different types of singularities and so the
computations would be very long and tedious. For our purpose it is
enough to have the following information on the ``shape'' of that number.

\begin{prop}\label{prop:numb}
  For any odd prime $p$ we have
\[\#X(\mathbb F_{p})=1+\sum_{i=1}^{4}\sum_{j=1}^{b_{2i}}\left(\frac {a_{i,j}}p\right)p^{i}+p^{5}-a_{p},\]
where the $a_{i,j}$ are square-free non-zero integers.
\end{prop}

\begin{proof}[Sketch of the proof of Proposition~\ref{prop:numb}]
  Using Ahlgren's result we only have to take into account the effect of adding the hypersurface
  at infinity and of all blow-ups. All these varieties 
  are resolutions of certain double covers branched along divisors of
  small degrees. Using the projection formula for finite maps, 
  it is not difficult to show that the Hodge spaces contributing to the odd
  cohomology groups are all zero. The even cohomology groups are spanned by 
  algebraic cycles, which project (under the
  double covering) onto cycles defined over $\mathbb Q$. 
  Consequently, the even cohomology groups can be generated by cycles which are either
  defined over $\mathbb Q$ or over some quadratic extension.
  Fix a non-invariant
  irreducible algebraic subvariety $Z$ and denote by $Z'$ its image under
  the involution defined by the double cover. 
  Clearly $Z+Z'$ is defined over $\QQ$. Assume that $Z$ (and hence also
  $Z'$) is defined over a quadratic extension $\QQ(\sqrt a)$. Recall that $p$ is an odd prime, 
  and hence a prime of good reduction.
  Over $\bar\FF_{p}$ the sum $(Z+Z')_p$ splits into a sum of two cycles $Z_{p}$ and
  $Z'_{p}$. Frobenius 
  maps the class $Z_{p}$ to the class of  $p^{i}Z_{p}$ or
  $p^{i}Z'_{p}$ ($i=5- \dim Z$) depending on whether $a$ is a square in
  $\FF_{p}$ or not. Consequently the class of the cycle $Z_{p}-Z'_{p}$ is an eigenvector with
  eigenvalue $(\frac ap)p^{i}$, where $(\frac ap)$ is the Legendre
  symbol. 
  Using the Lefschetz
  fixed-point formula we obtain the proposition.  
\end{proof}

\begin{proof}[Proof of Theorem \ref{theo:Ahlgren}] 
For every $1\le i\le4$ and $1\le j\le b_{2i}$ we consider the
one-dimensional Galois  representation $\rho_{i,j}$ with eigenvalues
$\left(\frac {a_{i,j}}p\right)p^{i}$ and define $\tilde \rho$ to be the
direct sum of all $\rho_{i,j}$. So $\tilde \rho$ is the Galois
representation associated to the algebraic cycles. 
Let $\bar\rho_{i}$ be the Galois action on the $i$-th cohomology 
and denote by $\bar\rho$ the direct sum of $\bar\rho_{2i}$,
$i=1,\dots,4$. 
Finally, denote by $\rho$ the Galois representation associated to the
unique cusp form of level 4 and weight 6. 
By Proposition~\ref{prop:numb} we can write the number of points of
$X(\FF_{p})$ as
\[1+p^{5}+\tr(\tilde\rho_{p})-\tr{\rho_{p}}.\]
By the Lefschetz fixed point
formula this is equal to
\[1+p^{5}+\tr(\bar\rho_{p})-\tr(\rho_{1,p})-\tr(\rho_{3,p})-\tr(\rho_{5,p})-\tr(\rho_{7,p})-\tr(\rho_{9,p}).\]
Comparing the above two formulas and clearing the signs we get 
\[\tr(\bar\rho_{p})+\tr(\rho_{p})=\tr(\tilde\rho_{p})+\tr(\rho_{1,p})+
\tr(\rho_{3,p})+\tr(\rho_{5,p})+\tr(\rho_{7,p})+\tr(\rho_{9,p})\]
So the representations \[\bar\rho\oplus\rho\] and
\[\tilde\rho\oplus\rho_{1}\oplus\rho_{3}\oplus\rho_{5}\oplus\rho_{7}\oplus\rho_{9}\]
have equal traces for any odd prime, and consequently 
they have isomorphic semi-sim\-pli\-fi\-cations (see \cite[Lemma p. I-11]{Serre}). Semi-sim\-pli\-fi\-ca\-tion
preserves the eigenvalues. By construction and  the Weil conjectures the
representation $\bar\rho\oplus\rho$ has no eigenvalue  with absolute
value equal to $p^{1/2}$ or $p^{3/2}$, and only two  eigenvalues  with absolute value equal
equal to $p^{5/2}$. So $H^{1}(X)=H^{3}(X)=H^{7}(X)=H^{9}(X)=0$ 
and the Galois representations $\rho$ and $\rho_{5}$ have
equal eigenvalues and hence isomorphic semi-sim\-pli\-fi\-ca\-tions.
\end{proof}

\begin{remark}
  The Ahlgren variety is birational to the quotient of the fourfold fiber
  product of the Legendre family, resp. the extremal rational elliptic surface with three singular fibers
  of Kodaira types $I_{2}, I_{2}, I^{*}_{2}$ (which in \cite
  {MirandaPersson} is denoted by $X_{222}$) by the group $\mz_{2}^{3}$.
  In each fiber this is the construction described in
  Section~\ref{sec:kummer} so  it is fibered by Calabi-Yau
  $4$-folds. 
\end{remark}

\subsection{Resolution of singularities of double arrangements}
\label{sec:arr}
In this subsection we shall describe in detail the procedure which we use
to resolve the singularities of Ahlgren's fivefold. 
Let $Y$ be an $n$-dimensional smooth projective manifold.
\begin{definition}\label{def:arr}
  A sum $D=\bigcup\limits_{i=1}^{N}D_{i}$ of smooth hypersurfaces
  $D_{i}$ in  $Y$ is called an \emph{arrangement} if for
  each subset $\{i_{i},\dots,i_{r}\}\subset\{1,\dots,N\}$ the
  (ideal-theoretic) intersection
  $C_{i_{1},\dots,i_{r}}=D_{i_{1}}\cap\dots\cap D_{i_{r}}$ is
  smooth. 
\end{definition}

The following lemma is obvious from the definitions.

\begin{lemma}
  Let $D=D_{1}\cup\dots\cup D_{N}\subset Y$ be an arrangement. Then
  \begin{enumerate}
  \item If  $\dim (D_{i_{1}}\cap\dots\cap D_{i_{r}})=n-r$ for some
    $\{i_{i},\dots,i_{r}\}\subset\{1,\dots,N\}$,  then
    $D_{i_{1}}\dots D_{i_{r}}$ intersect transversally.
  \item For any  $\{i_{i},\dots,i_{r}\}\subset\{1,\dots,N\}$ the tangent
    space to the intersection $D_{i_{1}}\cap\dots\cap D_{i_{r}}$ (at
    any point)
    equals the intersection of the tangent spaces to the divisors $D_{i}$.
  \end{enumerate}
\end{lemma}

We now consider the decomposition of the singular locus of $D$ by
multiplicities. For this we take the set $\mathcal S$ of all
components $C$ of 
intersections $D_{i_{1}} \cap \dots \cap D_{i_{r}}$ where $r \geq 2$ and
$\{i_{i},\dots,i_{r}\}\subset\{1,\dots,N\}$. To each element $C \in \mathcal S$ we assign its
multiplicity $m(C)=\mult _{C}D=\#\{i:C\subset D_{i}\}$ and
dimension $d(C)=\dim C$. An element $C\in \mathcal S$ will be
called \emph{near-pencil} if it is contained in an element
$C'\in \mathcal S$ with $d(C)=d(C')-1$ and $m(C)=m(C')+1$ (i.e $C$ is
cut-out from $C'$ by a single hypersurface).

If the arrangement $D\subset Y$ is even (as an element of the Picard
group $\Pic(Y)$), then there exists a double cover $\pi:X\lra Y$ of
$Y$ branched along $D$. Such a double cover is uniquely
determined by fixing a line bundle $\mathcal L$ on $Y$ 
with $\mathcal O(D)\mathcal\cong \mathcal L^{\otimes2}$. 
\begin{prop}\label{prop:arr}
Assume that for every singular variety $C\in\mathcal S$ either $C$ is
near-pencil or $\left\lfloor\frac
  {m(C)}2\right\rfloor=n-d(C)-1$ then $X$ admits a projective crepant
resolution of singularities. 
\end{prop}

\begin{proof}
Let $C\in \mathcal S$ be of dimension $d(C)=d$ and
multiplicity  $m(C)=m$. By the definition of an arrangement, this is a smooth
subvariety of $Y$ and we consider the blow-up 
\[\sigma:\tilde Y\lra Y\]
of $Y$ along $C$ with exceptional divisor $E$. Recall that $C$, and hence $E$, are irreducible by the
definition of $\mathcal S$. The pullback $\sigma^{*}D$ of  $D$ to $\tilde Y$ is
even in the Picard group of $\tilde Y$, but it is in general not
reduced. We define $D^{*}$ as the unique reduced and even divisor
satisfying \[\tilde D\leq D^{*}\leq \sigma^{*}D,\] where $\tilde D$ is
the strict transform of $D$. In fact
$D^{*}$ is equal to $\tilde D$ or $\tilde D+\varepsilon E$ where $\varepsilon = 0$ if $m$ is even and 
$\varepsilon=1$ if $m$ is odd. This means that, when the multiplicity is even we take
the strict transform of the branch locus as the new branch locus,
whereas when the multiplicity is odd we add the exceptional divisor.
Equivalently $D^{*}=\sigma^{*}D-2\lfloor\frac
  m2\rfloor E$. 
We have $K_{\tilde Y}+\frac12D^{*}=\sigma^{*}(K_{Y}+\frac12D)+(n-d-1-\left\lfloor\frac
  m2\right\rfloor)E$, and so $K_{\tilde
  Y}+\frac12D^{*}=\sigma^{*}(K_{Y}+\frac12D)$ exactly when $\left\lfloor\frac
  m2\right\rfloor=n-d-1$. We shall call a blow-up for which this
equality holds \emph{admissible}.

Assume now that $C\in\mathcal S$ is a minimal element (with respect to inclusion) among those 
components which are not
near pencil. Then, by assumption, the blow-up $\sigma$ along $C$ is admissible.  
We want to show that   $D^{*}$ is again an
arrangement satisfying the assumptions of the theorem. 
Let $D_{1},\dots,D_{k}$ be
the components of $D$ that contain $C$. Let us pick some other
components $D_{k+1},\dots,D_{k+p}$ and denote by $C_{1}$ the
intersection $C_{1}=C\cap D_{k+1}\cap\dots\cap D_{k+p}$. 
As the problem is local we can assume that $C_{1}$ is irreducible.  
Our aim is to
show that the intersection $\tilde D_{1}\cap\dots\cap \tilde D_{l}\cap
\tilde D_{k+1}\cap\dots\cap \tilde D_{k+p}$ is smooth. 

The intersection consists of two parts, namely the strict transform of the
intersection and
the intersection of the exceptional loci. 
The dimension of the former is 
less than or equal to $\dim C_{1}+\codim C-1-l$, and so its codimension
is greater than or equal to $\codim C_{1}-\codim C+1+l$. 
Since all the intersections of $C$ with $D_{k+j}$ are near pencil we
obtain that $\codim C_{1}-\codim C=p$ and that the codimension of the
intersection of the exceptional loci is greater than $p+l$, and hence this
is not a component of the intersection  $\tilde D_{1}\cap\dots\cap \tilde D_{l}\cap
\tilde D_{k+1}\cap\dots\cap \tilde D_{k+p}$. Consequently,  the intersection
$\tilde D_{1}\cap\dots\cap \tilde D_{l}\cap 
\tilde D_{k+1}\cap\dots\cap \tilde D_{k+p}$ equals the strict
transform of the intersection $ D_{1}\cap\dots\cap  D_{l}\cap
 D_{k+1}\cap\dots\cap  D_{k+p}$, and hence is smooth. To conclude that
 $D^{*}$ is an arrangement in the case of $m$ odd,  we also have to take
 the exceptional divisor of the blow-up into account. But this is
 transversal to any strict transform. 

To show that the arrangement $D^{*}$ satisfies the assumption of the
proposition, we observe that in the case of $m$ even the
exceptional varieties for $D^{*}$ are blow-ups of the exceptional
varieties for $D$, with the same multiplicities and dimensions. In the
case of $m$ odd, we have to add the intersections with the exceptional
divisors, but these are near-pencil singularities.

A resolution of singularities of $X$ can now be obtained by blowing-up all
the components $C\in\mathcal S$ which are not
near-pencil, starting from the smallest dimension. Since every blow-up
decreases the number of not near-pencil elements, the process
will terminate. As the intersection of two hyperplanes cannot be near-pencil,
the components of the final branch locus must be disjoint, and hence we get a
resolution of singularities.
Finally, since all blow-ups are admissible, the resulting
resolution is crepant. 

Denote by $\sigma:\tilde Y\lra Y$ the composition of all inverse maps to the blow-ups,
and by $\pi:X\lra Y$ (resp. $\tilde\pi:\tilde X\lra \tilde Y$) the
double cover of $Y$ (resp. of $\tilde Y$) branched along the divisor $D$
(resp. along $\tilde D$). Then there exists a unique map $\tilde\sigma:\tilde
X\lra X$ making the following diagram com\-mu\-ta\-tive
\[
\begin{CD}
  \tilde X @>\tilde\sigma>> X\\
  @V\tilde\pi VV @VV\pi V\\
\tilde Y @>\sigma>> Y
\end{CD}
\]

So the constructed crepant resolution of $X$ is given by a proper,
birational morphism.
\end{proof}

Clearly, the resolution is in general not unique, but depends
on the order of the blow-ups.

\bigskip
%\bigskip

\noindent
S\l awomir Cynk$^{*}$ and Klaus Hulek\\
Institut f\"ur Mathematik (C)\\
Universit\"at Hannover\\
Welfengarten 1, 30060 Hannover, Germany\\
{\tt cynk@math.uni-hannover.de, hulek@math.uni-hannover.de}

%\bigskip
%\bigskip
%\bigskip
\bigskip

\noindent {\small$^{*}$On leave from:}\\[1mm]
Instytut Matematyki\\ Uniwersytet Jagiello\'nski\\ 
ul. Reymonta 4, 30-059 Krak\'ow, Poland

\begin{thebibliography}{99}
\bibitem{AOP}A.~Ahlgren, K.~Ono, D.~Penniston, \emph{Zeta functions 
of an infinite family of $K3$ surfaces}\/.   
Amer. J. Math.  124  (2002),  no. 2, 353--368.

\bibitem{ahlgren}S.~Ahlgren, \emph{The points of a certain fivefold
    over finite fields and the twelfth power of the eta function}\/,
  Finite Fields and Their Applications 8, 18--33 (2002). 

\bibitem{Bor} C.~Borcea, \emph{Calabi-Yau threefolds and complex
  multiplication}\/,  Essays on mirror manifolds,  489--502,
  Internat. Press, Hong Kong, 1992. 

\bibitem{Bor2} C.~Borcea, \emph{$K3$ surfaces with involution and
    mirror pairs of Calabi-Yau manifolds}\/,  Mirror symmetry, II,
  717--743, AMS/IP Stud. Adv. Math., 1, Amer. Math. Soc., Providence,
  RI, 1997. 

\bibitem{CM} S.~Cynk and C.~Meyer,  
\emph{Geometry and arithmetic of certain double octic Calabi-Yau manifolds}\/, 
 Canad. Math. Bull.  48  (2005),  no. 2, 180--194.
 
\bibitem{DM} L.~Dieulefait and J.~Manoharmayum,  \emph{Modularity of rigid
  Calabi-Yau threefolds over $\mathbb Q$}.  Calabi-Yau varieties and
  mirror symmetry (Toronto, ON, 2001),  159--166, Fields
  Inst. Commun., 38, Amer. Math. Soc., Providence, RI, 2003.  
 
\bibitem{FM} J.-M.~Fontaine and B.~Mazur \emph{Geometric Galois
  representations}\/, Elliptic curves, modular forms, \& Fermat's last
  theorem (Hong Kong, 1993),  41--78, Ser. Number Theory, I,
  Internat. Press, Cambridge, MA, 1995. 

\bibitem{HV1} K.~Hulek and H.~Verrill, \emph{On modularity of rigid
    and nonrigid Calabi-Yau varieties associated to the root lattice
    $A_4$}\/, Nagoya Math. J.  179 (2005), 103--146.
 

\bibitem{HV2} K.~Hulek, H.~Verrill, \emph{On the modularity
    of Calabi-Yau threefolds containing elliptic ruled surfaces},
preprint (2005), math.AG/0502158, 
to appear in AMS/IP Studies in Advanced Mathematics ''Mirror Symmetry V'',
Proceedings of the BIRS Workshop on Calabi-Yau varieties
and mirror symmetry, December 6-11, 2003.


\bibitem{KS} 
H. H. Kim\ and\ F. Shahidi, Functorial products for ${\rm GL}\sb 2\times{\rm GL}\sb 3$ 
and the symmetric cube for ${\rm GL}\sb 2$. 
With an appendix by C. J. Bushnell and G. Henniart.  
Ann. of Math. (2) 155 (2002), 837--893 (2002).

 \bibitem{Livne} R.~Livn\'e, \emph{Motivic orthogonal two-dimensional
representations of ${\rm Gal}(\overline{\mathbb Q}/{\mathbb Q})$}\/,   Israel
  J. Math. 92 (1995), no. 1-3, 149--156. 

\bibitem{LY}
R.~Livn\'e and  N.~Yui,
\emph{The modularity of certain non-rigid Calabi-Yau threefolds}\/,
math.AG/0304497.
To appear:  J. Math. Kyoto Univ.

\bibitem{bookofMeyer}C.~Meyer, \emph{Modular Calabi-Yau threefolds},
Fields Institute Monograph \textbf{22} (2005), AMS.

\bibitem{MirandaPersson}R.~Miranda,  U.~Persson,  \emph{On extremal
    rational elliptic surfaces}\/,  Math. Z.  193  (1986),  no. 4,
  537--558.
  
\bibitem{Ribet} K.~Ribet, \emph{ Galois representations attached to eigenforms with Nebentypus}\/,
Modular functions of one variable, V (Proc. Second Internat. Conf.,
Univ. Bonn, Bonn, 1976), pp. 17--51. Lecture Notes in Math., Vol. 601,
Springer, Berlin, 1977.   

\bibitem{Serre} J.-P.~Serre, \emph{Abelian $l$-adic representations and
  elliptic curves}\/, With the collaboration of Willem Kuyk and John
  Labute. Revised reprint of the 1968 original. Research Notes in
  Mathematics, 7. A K Peters, Ltd., Wellesley, MA, 1998. 

\bibitem{SI}
T.~Shioda,  H.~Inose, 
\emph{On singular $K3$ surfaces}\/,
Complex analysis and algebraic geometry, pp. 119--136.
Iwanami Shoten, Tokyo, 1977. 

\bibitem{Ueno} K.~Ueno, \emph
{Classification theory of algebraic
    varieties and compact complex spaces}\/,  Lecture Notes in
  Mathematics, Vol. 439. Springer-Verlag, Berlin-New York, 1975. 


\bibitem{Voisin} C.~Voisin, \emph{Miroirs et involutions sur les
    surfaces $K3$}\/,   Journ\'ees de G\'eom\'etrie Alg\'ebrique d'Orsay (Orsay,
  1992).   Ast\'erisque  No. 218 (1993), 273--323.

\bibitem{Wiles} A.~Wiles, \emph{Modular elliptic curves and Fermat's
    last theorem\/}\/,   Ann. of Math.~(2)  141  (1995),  no. 3, 443--551.

\bibitem{Yui} N.~Yui, \emph{Update on the modularity of Calabi-Yau varieties}\/,
With an appendix by Helena Verrill. Fields Inst. Commun., 38, 
Calabi-Yau varieties and mirror symmetry (Toronto, ON, 2001), 307--362,
Amer. Math. Soc., Providence, RI, 2003.  

\end{thebibliography}
\end{document}